\documentclass{ws-jktr_ArXiv}
\usepackage{graphicx}
\begin{document}
\markboth{Sofia Lambropoulou, Stathis Antoniou}
{}

\catchline{}{}{}{}{}

\title{Topological Surgery, Dynamics\\ and Applications To Natural Processes}

\author{Sofia Lambropoulou} 

\address{Department of Mathematics, National Technical University of Athens, Athens, Greece. \\
sofia@math.ntua.gr}

\author{Stathis Antoniou}

\address{Department of Mathematics, National Technical University of Athens, Athens, Greece. \\
santoniou@math.ntua.gr}

\maketitle

\begin{abstract}
In this paper we observe that 2-dimensional 0-surgery occurs in natural processes, such as tornado formation and other phenomena reminiscent of hole drilling. Inspired by such phenomena, we introduce new theoretical concepts which enhance the formal definition of 2-dimensional 0-surgery with the observed dynamics. To do this, we first present a schematic model which extends the formal definition to a continuous process caused by local forces. Next, for modeling phenomena which do not happen on  surfaces but are three-dimensional, we fill the interior space by defining the notion of solid 2-dimensional 0-surgery. Finally, we connect these new theoretical concepts with a dynamical system and present it as a model for both 2-dimensional 0-surgery and natural phenomena exhibiting it. We hope that through this study, topology and dynamics of many natural phenomena will be better understood.
\end{abstract}


\ccode{2010 Mathematics Subject Classification: 57R65, 57N12, 57M25, 57M99, 37B99, 92B99.}

\section{Introduction}\label{Intro}

Topological surgery is a mathematical technique for creating new manifolds out of known ones. In \cite{1,2,3,4} we observed that surgery manifests itself in many natural processes of various scales. Examples of such phenomena comprise chromosomal crossover, magnetic reconnection, mitosis, gene transfer, the creation of Falaco Solitons, the formation of whirls and tornadoes, magnetic fields and the formation of black holes. In order to topologically model these phenomena, in \cite{4} we introduced the new notions of {\it dynamics} and {\it solid surgery} to the formal definition of topological surgery. In \cite{4}, we also connected these notions with a dynamical system by observing that its solutions are exhibiting 2-dimensional 0-surgery.

In this paper we focus on 2-dimensional 0-surgery and processes that resemble `hole drilling'. In Section {2}, we start by recalling the formal definition of surgery and we see how our new notions look like for this specific case. Then, in Section {3} we present the formation of tornadoes as an example of a natural process exhibiting this type of surgery. Finally, in Section {4} we present the aforementioned dynamical system which can serve as a model for both 2-dimensional 0-surgery and natural phenomena involving `hole drilling'. For a detailed treatment of natural processes exhibiting surgery in dimensions 1 and 2, for generalizations of the new notions of {\it dynamics} and {\it solid surgery} and for a connection between surgeries of different dimensions, the reader is referred to \cite{4}.

\section{Topological notions} \label{TopNo}
We start by recalling the formal definition of topological surgery and we see how our new notions enhance the definition of 2-dimensional 0-surgery.

\bigbreak 
{\noindent} {\bf The formal definition of surgery:} An \textit{m-dimensional n-surgery} is the topological procedure of creating a new $m$-manifold $M'$ out of a given $m$-manifold $M$ by removing a framed $n$-embedding $h:S^n\times D^{m-n}\hookrightarrow  M$, and replacing it with $D^{n+1}\times S^{m-n-1}$, using the `gluing' homeomorphism $h$ along the common boundary $S^n\times S^{m-n-1}$. That is:
\[M' = \chi(M) = \overline{M\setminus h(S^n\times D^{m-n})} \cup_{h|_{S^n\times S^{m-n-1}}} D^{n+1}\times S^{m-n-1}. \]
Note that from the definition we must have $n+1 \leq m$. For further reading, excellent references on the subject are \cite{5, 6, 7}.

\begin{figure}[ht!]
\begin{center}
\includegraphics[width=11cm]{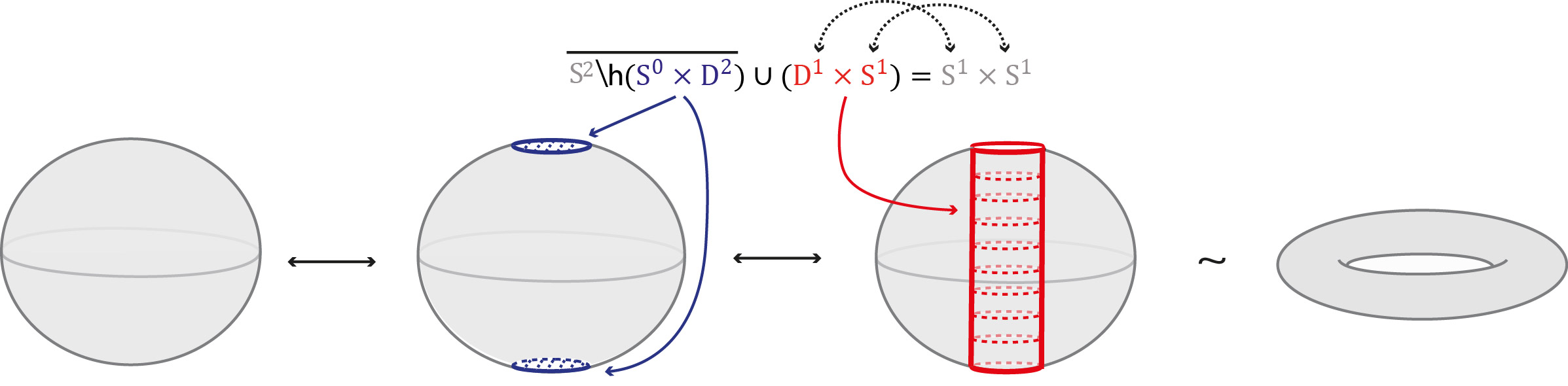}
\caption{Formal 2-dimensional 0-surgery.}
\label{2D0_Formal}
\end{center}
\end{figure}

For $m=2$ and $n=0$, starting with a 2-manifold $M$, we have the \textit{ 2-dimensional 0-surgery} whereby two discs $S^0\times D^2$ are removed from  $M$ and are replaced in the closure of the remaining manifold by a cylinder $D^1\times S^1$, which gets attached via a homeomorphism along the common boundary $S^0\times S^1$. The gluing homeomorphism of the common boundary can twist one or both copies of $S^1$. For $M=S^2$, the above operation changes the homeomorphism type from the 2-sphere to that of the torus, see Fig.~\ref{2D0_Formal}.

\smallbreak

\begin{figure}[ht!]
\begin{center}
\includegraphics[width=11cm]{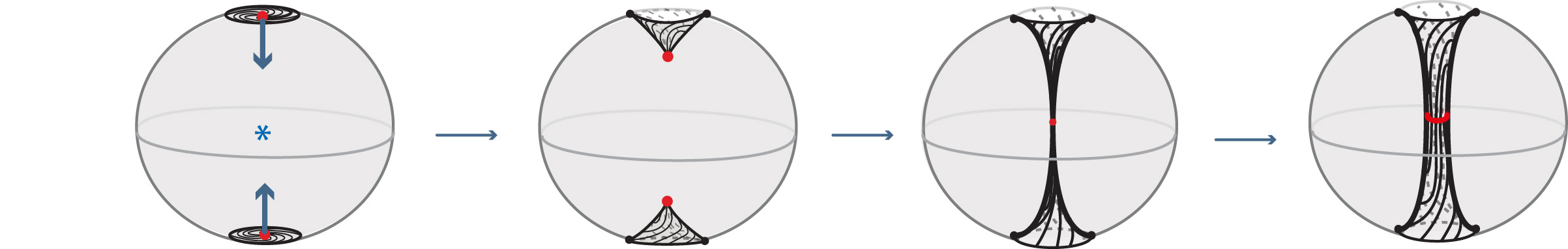}
\caption{Dynamic 2-dimensional 0-surgery.}
\label{2D0_Dynamics}
\end{center}
\end{figure}

{\noindent} {\bf Introducing dynamics:} In order to model topologically phenomena exhibiting 2-dimensional 0-surgery or to understand 2-dimensional 0-surgery through continuity we introduce dynamics to the above formal definition. As illustrated in Fig.~\ref{2D0_Dynamics}, in the case of $M=S^2$, we extend 2-dimensional 0-surgery to a continuous process caused by local forces and present the intermediate steps between the initial and the final manifold. 
 
\smallbreak
 
2-dimensional 0-surgery starts with two points, or poles, (in red) specified on the manifold  on which attracting forces created by an attracting center are applied (in blue). Then, the two discs $S^0\times D^2$, neighbourhoods of the two poles, approach each other. When the centers of the two discs touch, recoupling takes place and the discs get transformed into the final cylinder $D^1\times S^1$. The third instance of Fig.~\ref{2D0_Dynamics} illustrates the `recoupling', the singular stage at which topology changes. Note that these intermediate steps can be also explained through Morse theory, see \cite{8}.

It is worth pointing out that in Fig.~\ref{2D0_Formal} we used the standard embedding, while in Fig.~\ref{2D0_Dynamics}, the two discs $S^0\times D^2$  are embedded via a twisted homeomorphism $h_t$. More specifically, if we define the homeomorphisms  $\omega_1,\omega_2:D^{2}\to D^{2}$ to be rotations by $4\pi/3$ and $-4\pi/3$ respectively, then $h_t$ is defined as the composition $h_t:S^{0}\times D^2 \xrightarrow{\omega_1 \amalg  \omega_2}  S^0\times D^{2}  \xrightarrow{h} M $. We used an example of a twisted homeomorphism $h_t$ as it is more commonly found in natural processes.
 
\bigbreak 
{\noindent} {\bf Solid surgery:} In order to model natural phenomena undergoing 2-dimensional 0-surgery which do not happen on surfaces but are three-dimensional, we introduce the notion of solid 2-dimensional 0-surgery. The interior of the initial manifold is now filled in and the interior space undergoes the same type of surgery as the initial manifold. 

\smallbreak
\begin{figure}[ht!]
\begin{center}
\includegraphics[width=13cm]{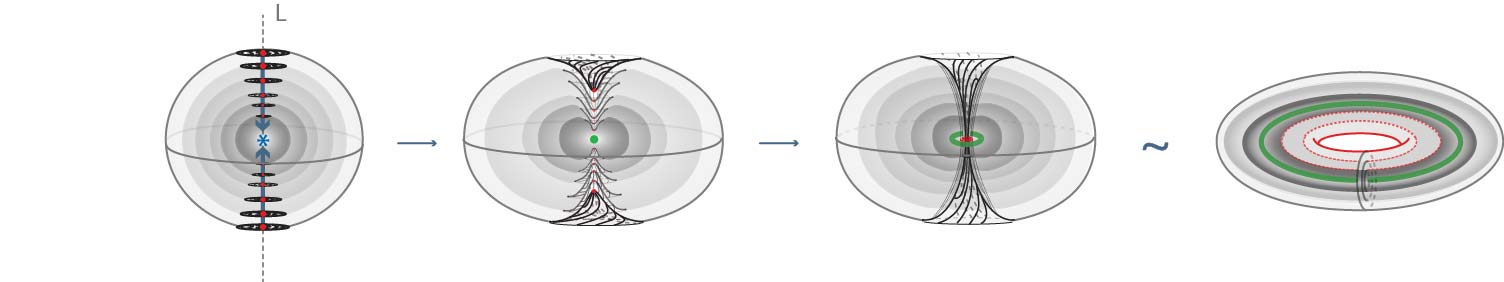}
\caption{Solid 2-dimensional 0-surgery.}
\label{2D0_Solid}
\end{center}
\end{figure}

More specifically, we start with the $3$-ball of radius 1 with polar layering: 
$$
D^3 = \cup_{0<r\leq 1} S^2_r \cup \{P\},
$$ 
where $r$ the radius of the 2-sphere $S^2_r$ and $P$ the limit point of the spheres, that is, their common center and the center of the ball. \textit{Solid 2-dimensional 0-surgery on $D^3$} is the  topological procedure whereby 2-dimensional 0-surgery takes place on each spherical layer that $D^3$ is made of. More precisely, as illustrated in Fig.~\ref{2D0_Solid}, on all spheres $S^2_r$  colinear pairs of antipodal points are specified, on which the same colinear attracting forces act. The poles have disc neighborhoods of analogous areas. Then, 2-dimensional 0-surgeries are performed on the whole continuum of the concentric spheres using the same homeomorphism $h$. 

Moreover, 2-dimensional 0-surgery on the limit point $P$ is defined to be the limit circle of the nested tori resulting from the continuum of 2-dimensional 0-surgeries. That is, the effect of \textit{2-dimensional 0-surgery on a point is defined to be the creation of a circle}.

The process is characterized on one hand by the 1-dimensional core $L$ of the solid cylinder which joins the two selected antipodal points of the outer shell and intersects each spherical layer at its two corresponding antipodal points, and on the other hand by the homeomorphism $h$. The process results in a continuum of layered tori and can be viewed as drilling out a tunnel along $L$ according to $h$. For a twisted homeomorphism $h_t$, which is the case shown in Fig.~\ref{2D0_Solid}, this agrees with our intuition that, for opening a hole, \textit{drilling with twisting} seems to be the easiest way.


\section{Tornadoes and `hole drilling' processes} \label{Tornadoes}

\begin{figure}[ht!]
\begin{center}
\includegraphics[width=12cm]{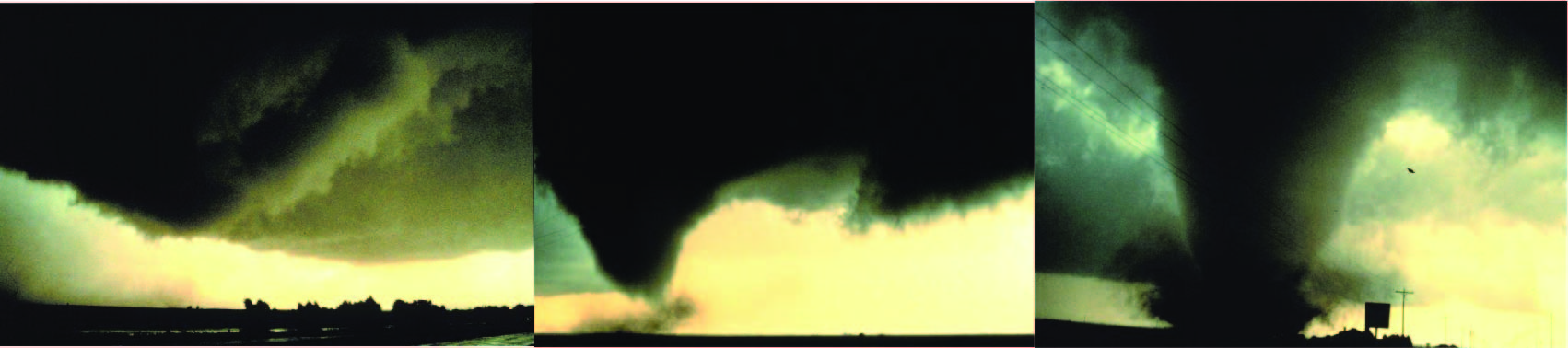}
\caption{A sequence of images showing the birth of a tornado.}
\label{TornadoSeq}
\end{center}
\end{figure}

As mentioned in the Introduction, the dynamics and the theoretical forces introduced in our schematic topological models are also observed in a multitude of natural processes. Consider for example the similarities between the dynamics shown in the topological schematic model of Fig.~\ref{2D0_Solid} and the formation of a tornado in Fig.~\ref{TornadoSeq}. 

\smallbreak

Indeed, if certain meteorological conditions are met, an attracting force between the cloud and the earth beneath is created. First, the rotating cloud base lowers. This lowering becomes a funnel-shaped cloud which starts descending toward the ground. Finally, the visible funnel extends to the ground. In analogy to 2-dimensional 0-surgery, first the poles are chosen, one on the tip of the cloud and the other on the ground, and they seem to be joined through an invisible line. Then, starting from the first point, the wind revolves in a helicoidal motion toward the second point, resembling `hole drilling' along the line, until the hole is drilled. Note that, in this case we have a realization of solid 2-dimensional 0-surgery, where the attracting center coincides with the ground and we only see helicoidal motion in one direction.  
 
As the same process happens on all layers of wind, the process of tornado formation is a solid 2-dimensional 0-surgery. As each layer revolves in a helicoidal motion, the two discs $S^0\times D^2$ are always embedded via a twisted embedding $h_t$. The hole is drilled along the line $L$ shown in Fig.~\ref{2D0_Solid}. Once the tornado is formed, the line $L$ coincides with the tornado's core line.

The initial manifold can be considered as being a 3-ball surrounding the phenomenon. The physical meaning of this consideration is that the process is triggered by the difference in the conditions of the lower and upper atmosphere, hence the initial manifold is considered as the 3-ball containing this air cycle. 

Further, we note that the dynamics of our general schematic topological model of Fig.~\ref{2D0_Solid} are also found in whirlwinds, in the creation of Falaco Solitons (see \cite{9}  
 for details) in the creation of vortices, and can be more generally applied to natural processes exhibiting a `hole drilling' behavior.

\section{A dynamical system for 2-dimensional 0-surgery} \label{DS} 
We shall now discuss the following three-dimensional dynamical system presented in \cite{10}:

\begin{equation} 
\left\{
\begin{array}{l}
\frac{dX}{dt}=X-XY+CX^2-AZX^2 \\\\

\frac{dY}{dt}=-Y+XY \\\\

\frac{dZ}{dt}=-BZ+AZX^2 \\
\end{array}
\right\} \ A, B, C > 0
\end{equation} 
\bigbreak
 
This system is a two-predator and one-prey generalized Lotka-Volterra model, where the predators $Y,Z$ do not interact directly with one another but compete for prey~$X$. In \cite{10}, the authors analyze parameters $A,B,C$ in order to determine the bifurcation properties of the system, that is, to determine the changes in the qualitative or topological structure of the solutions of this family of differential equations. Conditions for which the ecosystem of the three species results in steady, periodic or chaotic behavior are determined. As it turns out, the various stability conditions can be determined by only two parameters: $C$ and $B/A$. It is also shown in \cite{10} that stable solutions are generated left of and including the line $B/A=1$ while chaotic/periodic regions appear on the right of the line $B/A=1$. 

\smallbreak
\begin{figure}[ht!]
\begin{center}
\includegraphics[width=12cm]{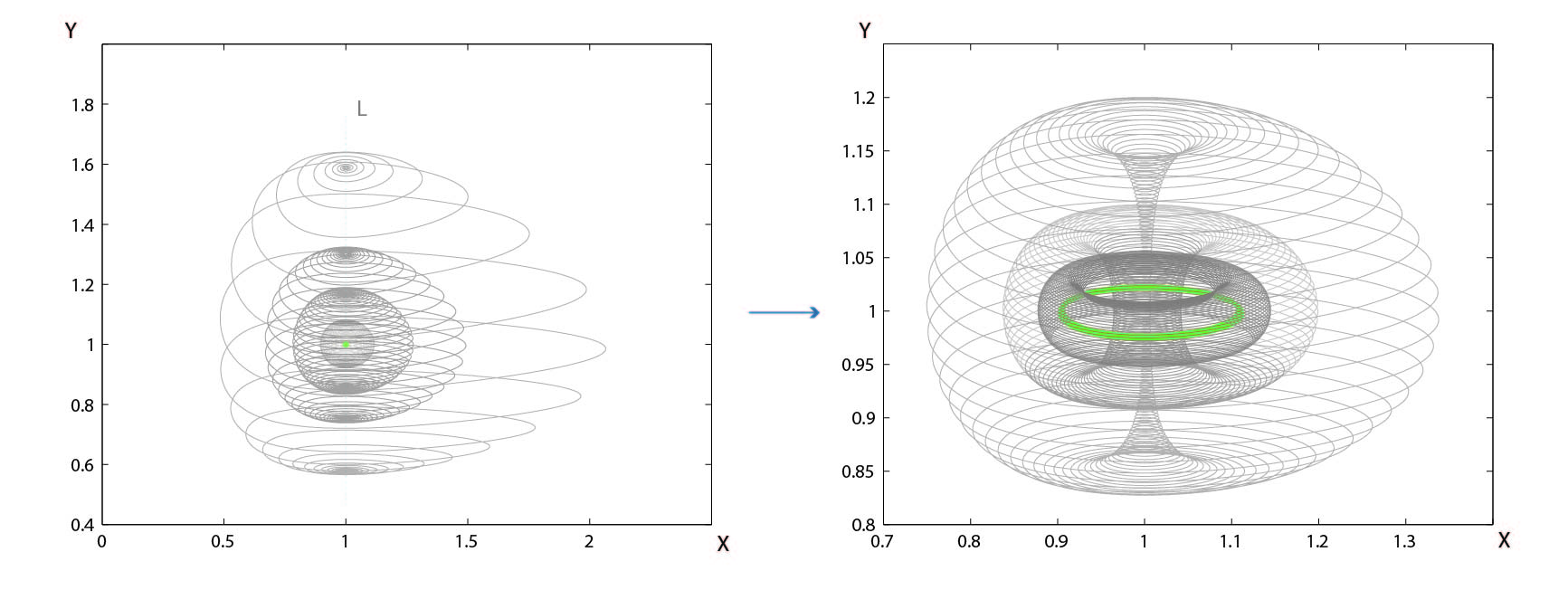}
\caption{The trajectories of the dynamical system are performing solid 2-dimensional 0-surgery.}
\label{LV_Surgery_jktr}
\end{center}
\end{figure}

In subsequent work \cite{11} the authors show that the solutions of this system form a nesting of concentric spheres in the stable region of $B/A=1$. As $B/A$ becomes greater than $1$,  spherical trajectories follow a motion reminiscent of `hole drilling' along a segment $L$ and turn into toroidal trajectories. Segment $L$ connects all spherical layers and the authors refer to it in \cite{10} as the `slow manifold' due to the fact that trajectories move slower when passing near it. In Fig.~\ref{LV_Surgery_jktr}, we reproduce the numerical simulations done in \cite{11}. More precisely, the spherical nesting is formed by trajectories  initiated at points $[1,1.59,0.81]$, $[1,1.3,0.89]$, $[1,1.18,0.95]$, $[1,1.08,0.98]$ and $[1,1,1]$ for parameters $A=B=C=3$, while the toroidal nesting is manifested by trajectories initiated at points $[1.1075,1,1]$, $[1,1,0.95]$, $[1,1,0.9]$ and $[1,1,1]$ for $A=2.9851$ and $B=C=3$. Comparing with Fig.~\ref{2D0_Solid}, we establish that solid 2-dimensional 0-surgery is performed by changing the parameter space from $B/A=1$ to $B/A>1$. 

However, as the authors elaborate in \cite{10}, while for $B/A=1$ the entire positive space is filled with nested spheres, when $B/A>1$, only spheres up to a certain volume become tori. More specifically, quoting the authors: ``to preserve uniqueness of solutions, the connections through the slow manifold $L$ are made in a way that higher volume shells require slower, or higher resolution, trajectories within the bundle". As they further explain: ``higher volume shells rapidly collapse or dissipate". As shown in both \cite{10} and \cite{11}, the outermost shell of the toroidal nesting is a fractal torus. Note that in Fig.~\ref{LV_Surgery_jktr}, we do not show the fractal torus because we are interested in the interior of the fractal torus which supports a topology stratified with toroidal surfaces. 

\smallbreak
\begin{figure}[ht!]
\begin{center}
\includegraphics[width=12cm]{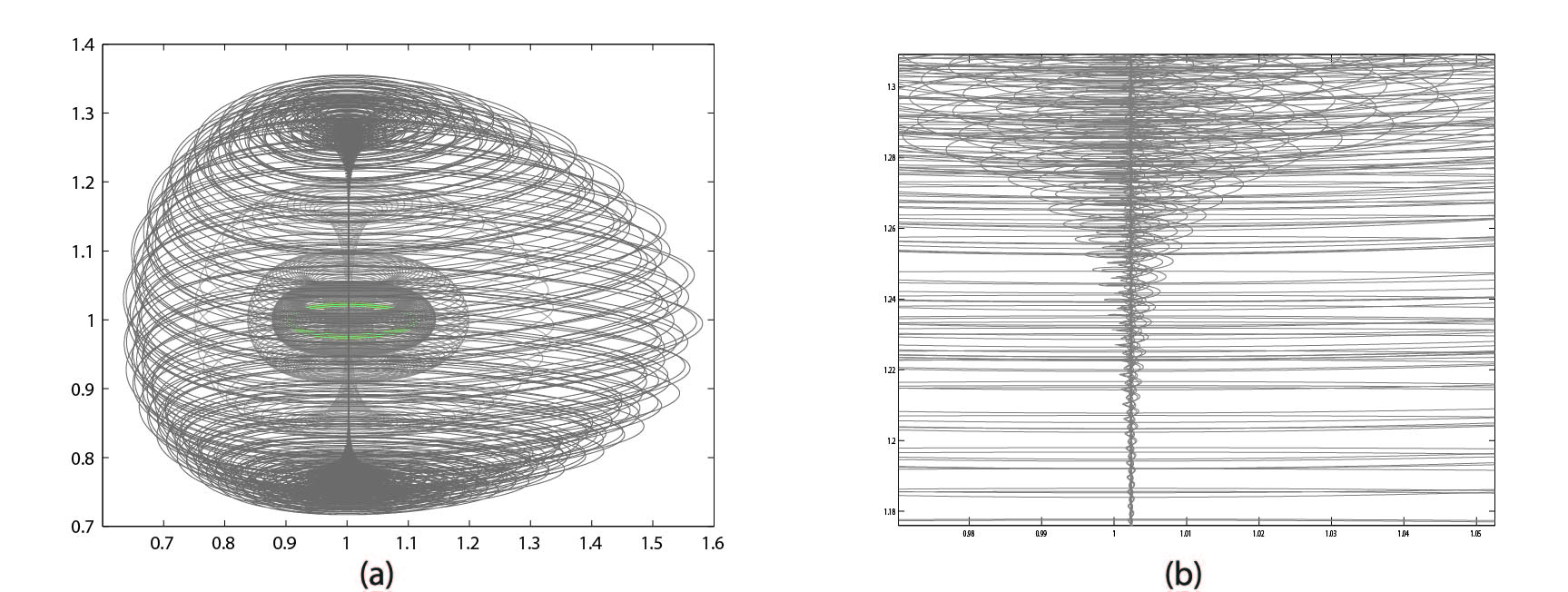}
\caption{(a) The fractal torus for $B/A>1$   (b) Zooming on the slow manifold of the fractal torus}
\label{LV_Surgery_jktr_2}
\end{center}
\end{figure}

In Fig.~\ref{LV_Surgery_jktr_2} (a), we present the outermost fractal torus as a trajectory initiated at point $[1.45,1,1.45]$ for the same parameters used in Fig.~\ref{LV_Surgery_jktr}, namely $A=2.9851$ and $B=C=3$. On the inside of the fractal torus, one can still see the periodic toroidal nesting. By zooming on the slow manifold of the outermost fractal torus shell, in Fig. 6 (b) we can view the ‘hole drilling’ behavior of the trajectories. Here, the solutions' indetermination creates trajectories that are more life-like than the periodic solutions of Fig.~\ref{LV_Surgery_jktr} and resemble even more natural processes such as the formation of tornadoes. Indeed, compare Fig.~\ref{LV_Surgery_jktr_2} (b) with Fig.~\ref{TornadoSeq}.

Therefore, except for being a model for solid 2-dimensional 0-surgery, this dynamical system constitutes a model for tornado formation and other natural processes exhibiting solid 2-dimensional 0-surgery and `hole drilling' behaviors.

\section{Conclusion} \label{Conclusion} 
In this paper, we first enhanced the static description of 2-dimensional 0-surgery by introducing dynamics, by means of attracting forces. We then filled in the interior space by introducing the notion of solid 2-dimensional 0-surgery and the definition of surgery on a point. We presented tornado formation as an example of natural process exhibiting these dynamics. Further, we connected solid 2-dimensional 0-surgery with a dynamical system. This connection gives us on the one hand a mathematical model for 2-dimensional 0-surgery and, on the other hand, a dynamical system modeling natural phenomena exhibiting 2-dimensional 0-surgery through a `hole-drilling' process.
 
We hope that through this study, topology and dynamics of similar natural phenomena, as well as topological surgery itself, will be better understood and that our connections will serve as ground for more insightful observations.

\section*{Acknowledgments}
We are grateful to Louis Kauffman and Cameron Gordon for many fruitful conversations on Morse theory and 3-dimensional surgery. We would also like to thank Nikola Samardzija for discussing the topological aspects of the dynamical system presented in \cite{10}. Finally, we would like to acknowledge a comment by Tim Cochran pointing out the connection of our new notions with Morse theory. 

This research has been co-financed by the European Union (European Social Fund - ESF) and Greek national funds through the Operational Program "Education and Lifelong Learning" of the National Strategic Reference Framework (NSRF) - Research Funding Program: THALES: Reinforcement of the interdisciplinary and/or inter-institutional research and innovation.

\end{document}